\documentclass[12pt]{amsart}
\usepackage{color}
\usepackage{epsfig}
\usepackage{verbatim}
\usepackage{clrscode3e}

\newtheorem{theorem}{Theorem}[section]

\newtheorem{lemma}[theorem]{Lemma}

\begin{document}
\title {A characterization of 2-neighborhood degree list of diameter 2 graphs}
\date{}
\maketitle 
\small
\begin{center}
\begin{tabular}{lll}
\shortstack{\normalsize
N. Benakli\\ 
Department of Mathematics \\  
NYCCT, CUNY\\
Brooklyn, NY 11201\\
nbenakli@citytech.cuny.edu 
}
&
\shortstack{\normalsize
E. Halleck  \\ 
Department of Mathematics \\  
NYCCT, CUNY\\
Brooklyn, NY 11201\\
ehalleck@citytech.cuny.edu
}
\end{tabular}
\end{center}
\begin{center}
\begin{tabular}{lll}
\shortstack{\normalsize
S. R. Kingan \\     
Department of Mathematics \\
Brooklyn College, CUNY\\
Brooklyn, NY 11210\\
skingan@brooklyn.cuny.edu
}
\end{tabular}
\end{center}

 \begin{abstract}
\normalsize
 Let $N_2DL(v)$ denote the set of degrees of vertices at distance 2 from $v$. The  $2$-neighborhood degree list of a graph is a listing of $N_2DL(v)$ for every vertex $v$. A degree restricted $2$-switch on edges $v_1v_2$ and $w_1w_2$, where $deg(v_1)=deg(w_1)$ and $deg(v_2)=deg(w_2)$, is the replacement of a pair of edges $v_1v_2$ and $w_1w_2$ by the edges $v_1w_2$ and $v_2w_1$ given that $v_1w_2$ and $v_2w_1$ did not appear in the graph originally.  Let $G$ and $H$ be two graphs of diameter 2  on the same vertex set. We prove that  $G$ and $H$ have the same  $2$-neighborhood degree list if and only if  $G$ can be transformed into $H$ by a sequence of degree restricted $2$-switches. 
\end{abstract}

\normalsize

\section {\bf Introduction}

Two graphs $G_1$ and $G_2$ are {\it isomorphic} if there is a bijection from $V(G_1)$ to $V(G_2)$ that preserves adjacencies. 
Two labeled graphs $G_1$ and $G_2$ are {\it label isomorphic} or {\it identical} if there is a bijection from $V(G_1)$ to $V(G_2)$ that preserves labeled adjacencies. The {\it degree} of a vertex $v$, denoted by $deg(v)$, is the number of edges incident to $v$.  
The {\it distance} between a pair of vertices $u$ and $v$ in a graph $G$, denoted by $d(u, v)$, is the length of the shortest path between $u$ and $v$. The {\it diameter}, denoted by diam(G),  is the maximum value of $d(u, v)$, where the maximum is taken over all pairs of vertices $u$ and $v$ in $G$.  The eccentricity of vertex $v$, denoted by $e(v)$, is the maximum value of 
$d(u, v)$, where the maximum is taken over all vertices $u \neq v$.

\bigskip
 Let $G$ be a graph with vertices $v_1, \dots , v_n$.  The {\it degree sequence} of a graph is a listing of the degrees of its vertices: $deg(v_1), deg(v_2), \dots , deg(v_n)$.
The set of all vertices adjacent to a vertex $v$ is denoted by $N(v)$ and called the {\it neighborhood} of $v$. Note that the neighborhood of $v$ does not contain $v$ itself.  For each  vertex $v$, the {\it neighborhood degree list} of $v$, denoted by $NDL(v)$, is the list of  degrees of vertices in $N(v)$. The {\it neighborhood degree list} of $G$, denoted by $NDL(G)$,  is  the list of lists $$\{NDL(v_1), NDL(v_2), \dots NDL(v_n)\}$$ 
By convention the degree sequence and the neighborhood degree list of a vertex are written in descending order. The concept of neighborhood degree list was introduced independently by Barrus and Donavan {\bf [\ref{BarrusDonovan2018}]} and Bassler {\it et al} {\bf [\ref{Bassleretal2015}]}. 

 \bigskip

In this paper we generalize the notion of neighborhood degree list.     Let $N_k(v)$ be the set of vertices of distance $k\ge 1$ from $v$. In this notation $N(v)=N_1(v)$.  Observe that $k\le diam(G)$. The {\it $k$-neighborhood degree list} of $v$, denoted by $N_kDL(v)$, is the list of degrees of vertices in $N_k(v)$. The {\it $k$-neighborhood degree list} of  $G$, denoted by $N_kDL(G)$,  is  the list of lists $$\{N_kDL(v_1), N_kDL(v_2), \dots N_kDL(v_n)\}.$$  Essentially we are considering concentric balls of vertices of increasing distance centered around a vertex.  

 \bigskip

The concept of $N_kDL$ is a strengthening of the well known distance degree sequence. For a vertex $v$, let $deg_0(v)=1$ and for $k\ge 1$, let $deg_k(v)=|N_k(v)|$. The {\it  distance degree sequence} of $v$  is the sequence
$$(deg_o(v), deg_1(v), deg_2(v), \dots , deg_{e(v)}(v)).$$
It was introduced by Randic in {\bf [\ref{Randic1979}]} to distinguish chemical isomers by their graph structure {\bf[\ref{BuckleyHarary1990}]}. A graph in which all the vertices have the same distance degree sequence is called {\it distance degree regular.}  Distance regular graphs are regular graphs. However, not all regular graphs are distance regular. Distance regular graphs have applications in numerous areas including algebraic combinatorics and coding theory. See {\bf [\ref{Godsil1995}]} and {\bf  [\ref{DamKoolen2016}]}. 

\bigskip
The motivation for our definition of $k$-neighborhood degree list comes from the problem of identifying fake followers on Instagram and Twitter. A New York Times article
titled ``The follower factory''\footnote {https://www.nytimes.com/interactive/2018/01/27/technology/social-media-bots.html}  explains how celebrities purchase followers from companies that create millions of such accounts and sell them as followers. Such companies are called follower factories. An Instagram influencer's follower count (i.e. degree in the social network) may be high, but the followers may be mostly vertices of low degree. Such fake influencers would be flagged by computing their $NDL$. In case the fake accounts have a degree greater than 1 in an effort to hide that they are fake accounts, then computing $N_kDL$, for $k\ge 1$, would reveal an anomaly in the pattern of $N_kDL$ lists. Measures of centrality like betweenness centrality, eigenvalue centrality, PageRank, etc. can also be used to compare vertices, but they are global measures designed for specialized applications. On the other hand $N_kDL$ is a local measure. In many cases the entire graph is unknown and a local measure of influence is needed.
\bigskip

Although the terminology is quite different, the notion of $N_kDL(v)$ appears in Roberio {\it et al} {\bf [\ref{Riberioetal2017}]}. 
Approaches for inferences on graphs rely on finding ways to embed vertices into the $n$-dimensional real vector space ($R^n$) so that ``similar'' vertices are embedded near each other. The approach described in {\bf [\ref{Riberioetal2017}]} uses the number of links of a vertex to its neighbors, number of links of the neighbors to their neigbors, and so on. In other words they use $N_kDL(v)$ as a measure of similarity between vertices and they establish experimentally that $N_kDL(v)$ is better than state-of-the-art techniques in capturing similarity of vertices.

\bigskip

In this paper we give a short proof of Barrus and Donovan's theorem characterizing graphs with the same $NDL$ and a characterization of diameter 2 graphs with the same $N_2DL$. 

\begin{figure}[h]
\centering
\includegraphics[width=3in]{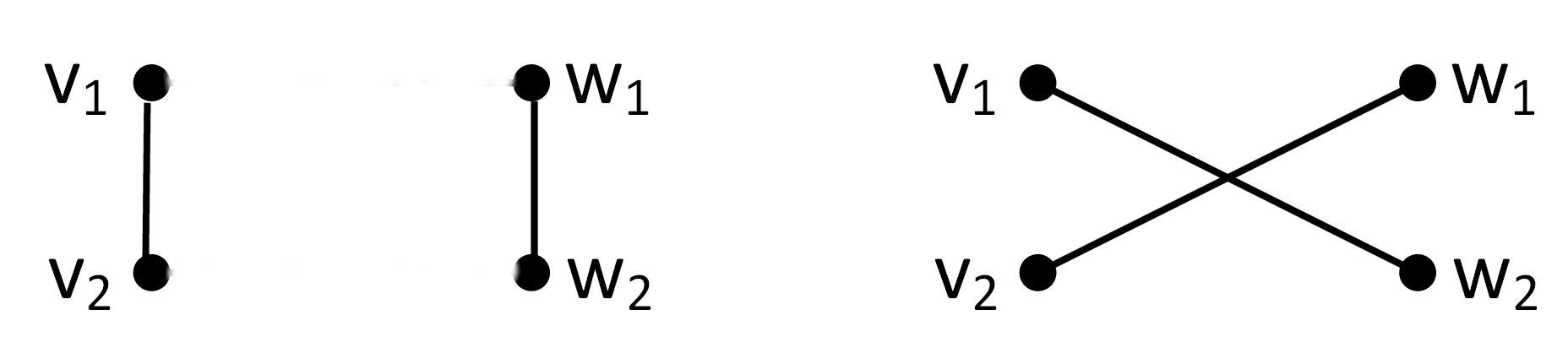}
\caption{The 2-switch operation 
\label{fig01}}
\end{figure}

\bigskip
A {\it  $2$-switch} in a graph  is the replacement of a pair of edges $v_1v_2$ and $w_1w_2$ by the edges $v_1w_2$ and $v_2w_1$ given that $v_1w_2$ and $v_2w_1$ did not appear in the graph originally. There may or may not be edges between pairs of vertices $v_1, w_1$  and $v_2, w_2$. The 2-switch operation is illustrated in  Figure \ref{fig01}.  It was the key idea in a result by Havel {\bf [\ref{Havel1955}]} and independently by Hakimi {\bf [\ref{Hakimi1962}]} that characterizes precisely the sequence of numbers that correspond to the degree sequence of a graph. See also {\bf [\ref{FulkersonHoffmanMcAndrew1965}]}. Observe that a 2-switch on a pair of edges does not alter the degrees of the four vertices involved. Therefore a 2-switch does not alter the degree sequence of the resulting graph. If some 2-switch turns $G$ into $G'$, then a 2-switch on the same four vertices turns $G'$ into $G$. 
\bigskip

\bigskip

A {\it degree restricted $2$-switch} on edges $v_1v_2$ and $w_1w_2$ is a 2-switch performed when $deg(v_1)=deg(w_1)$ and $deg(v_2)=deg(w_2)$.  (Barrus and Donavan call this an $n$-switch.) 
The next theorem is the main result in this paper.

 \bigskip
\begin{theorem} \label{diam2} Let $G$ and $H$ be two graphs on $n$ vertices with  diameter $2$. Then  $G$ and $H$ have the same  $2$-neighborhood degree list if and only if   $G$ can be transformed into $H$ by a sequence of degree restricted $2$-switches.
\end{theorem} 

 \bigskip


\section {The proof of Theorem 1.1}

 \bigskip

\begin{figure}[h]
\centering
\includegraphics[width=3in]{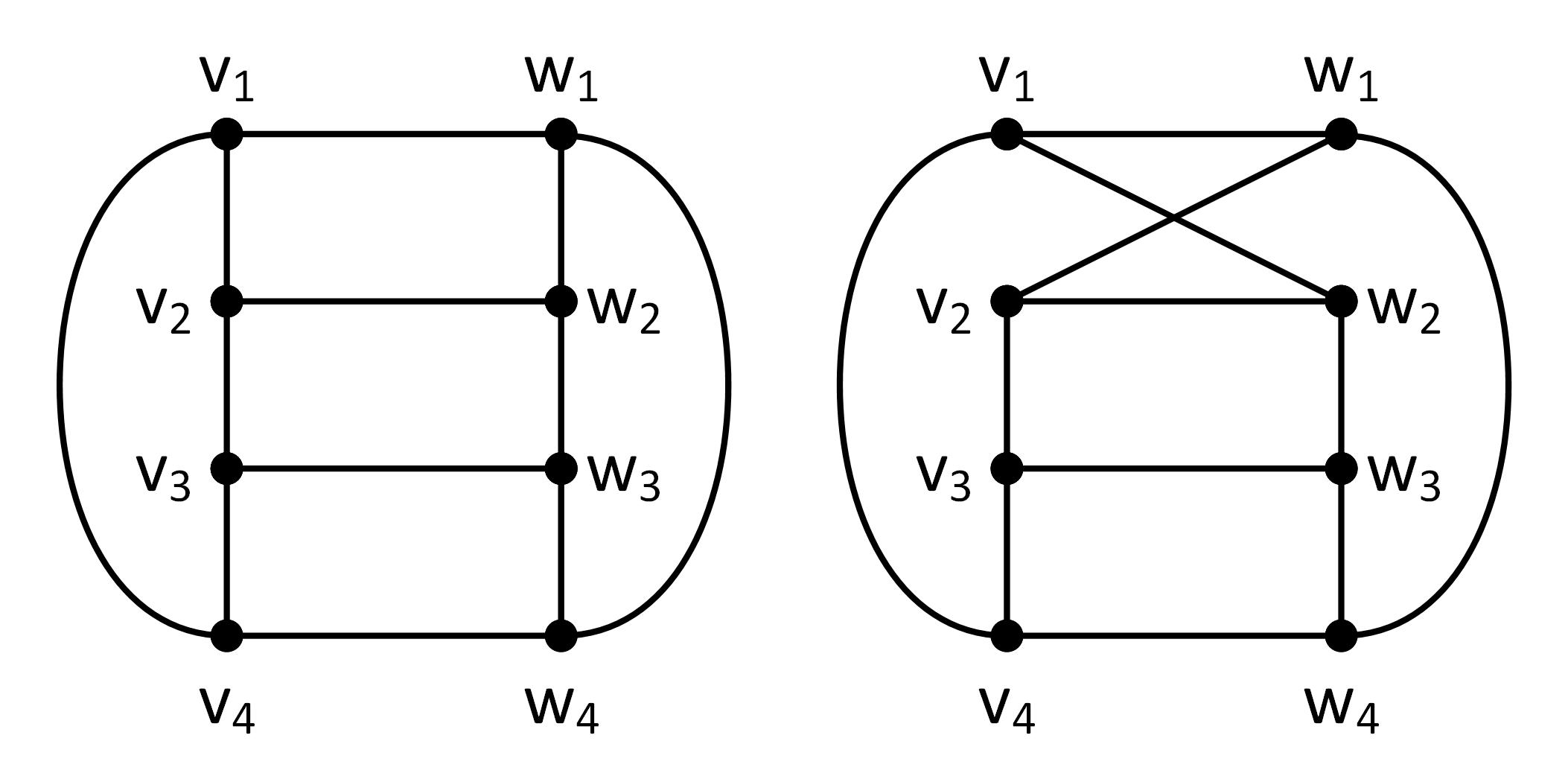}  
\caption{Circular $4$-ladder and Mobius $4$-ladder
\label{cube} }
\end{figure}

When generalizing $NDL$ to $N_2DL$ one problem that comes up  is that the 2-switch operation can alter the diameter. Consider for example a 2-switch performed on the cube graph (circular 4-ladder) that converts it to the Mobius 4-ladder as shown in Figure \ref{cube}. Observe that $NDL$ is preserved in both graphs. However, $N_2DL$ is not preserved. The diameter of the cube is 3, but when a 2-switch operation is done to obtain the Mobius 4-ladder, the diameter is reduced to 2. Thus the 2 graphs have different $N_2DL$.   

 \bigskip

Let $G$ and $H$ be two graphs on $n$ vertices. Berge proved that $G$ and $H$ have the same  degree sequence if and only if  $G$ can be transformed into $H$ by a sequence of $2$-switches.  See {\bf [\ref{West1995}, p. 47]}. We do not use this result, rather we use the technique that Berge uses. The next result appears in {\bf [\ref{BarrusDonovan2018}, Theorem 3.3]}. We give a different proof based on the lexicographic ordering of the neighborhood degree list. 
\bigskip

Let $a$ and $b$ be two vertices in a graph such that $deg(a)=deg(b)=t$. Let $N(a)=\{a_1, \dots , a_t\}$ and $N(b)=\{b_1, \dots , b_t\}$, where the vertices are listed in descending order based on their degrees.
 We say $NDL(a)=NDL(b)$ if the ordered lists of degrees are the same. In other words
$$(deg(a_1), \dots , deg(a_t)) =(deg(b_1), \dots , deg(b_t))$$  
 We say $NDL(a) < NDL(b)$ if using the lexicographic ordering $$(deg(a_1), \dots , deg(a_t)) <(deg(b_1), \dots , deg(b_t)).$$   Lexicographic ordering  is defined recursively. 
If $deg(a_1) < deg(b_1)$, then $NDL(a) < NDL(b)$. If $deg(a_1) = deg(b_1)$, then the  order is determined by the lexicographic order of $(deg(a_2), \dots , deg(a_t))$ and $(deg(b_2), \dots , deg(b_t))$. If $deg(a_2) < deg(b_2)$, then $NDL(a) < NDL(b)$. If $deg(a_2) = deg(b_2)$, then check the sequences $(deg(a_3), \dots deg(a_t))$ and $(deg(b_3), \dots deg(b_t))$, and so on.

\bigskip

\begin{lemma} \label{NDL} {\bf (Barrus and Donavan 2018)} Let $G$ and $H$ be two graphs on $n$ vertices. Then $G$ and $H$ have the same  neighborhood degree list if and only if $G$ can be transformed into $H$ by a sequence of degree restricted $2$-switches.
\end{lemma} 

\noindent {\bf Proof.} One direction is straightforward. If $G$ can be transformed into $H$ by a sequence of degree restricted 2-switches,
then clearly  $G$ and $H$ have the same neighborhood degree list.

 \bigskip
Conversely, suppose $G$ and $H$ have the same neighborhood degree list. The proof is by induction on $n\ge 4$. The result holds for graphs on 4 vertices trivially. Assume that the result holds for all graphs with $n-1$ vertices.

 \bigskip

Let $w$ be a vertex of maximum degree $\Delta$ in $G$. Let $z$ be a neighbor of $w$ and let $S$ be the set of all vertices that are not neighbors of $w$, but have the same degree as $z$. If $S=\phi$, then proceed to the next neighbor. Otherwise suppose $S\neq \phi$. Choose $x\in S$ so that $NDL(x)$ is highest among vertices of $S$. If $NDL(x)\le NDL(z)$, then again proceed to the next neighbor of $S$.

\begin{figure}[h]
\centering
\includegraphics[width=1.5in]{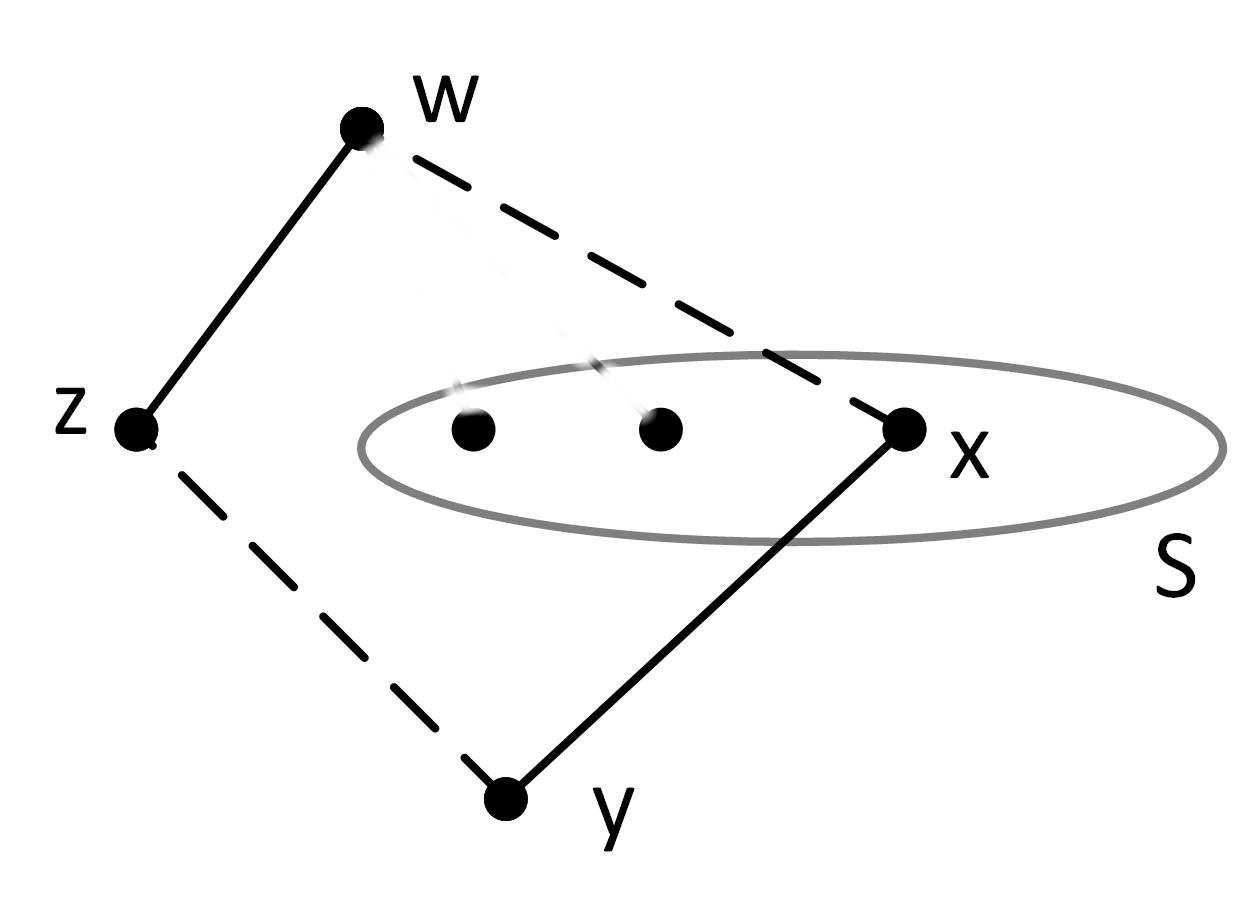}
\caption{Degree restricted 2-switch
\label{fig03}} 
\end{figure}

 \bigskip
Suppose $NDL(x) > NDL(z)$. Note that $w$ is a vertex of maximum degree. Since $w$ is  incident to $z$, but not to $x$, there exists $y$ incident to $x$, but not to $z$, such that $deg(y)=deg(w)$. Thus we can perform a degree restricted 2-switch operation on $wz$ and $yx$. Delete $wz$ and $yx$ and add $wx$ and $yz$. (See Figure \ref{fig03}.)  Observe that the resulting graph has the same $NDL$ as $G$ (and consequently the same degree sequence).

 \bigskip

Repeat the above process for every neighbor of $w$  to obtain a graph $G^*$  where $NDL(G^*)=NDL(G)$  and neighbors of $w$ in $G^*$ are adjacent to vertices with the same degrees as neighbors of $w$ in $G$, but with highest $NDL$.

 \bigskip
Similarly, choose a vertex $w_H$ in $H$ of highest degree such that $NDL_H(w_H)=NDL_G(w)$. Such a vertex exists since $G$ and $H$ have the same $NDL$. There exists a sequence of degree restricted 2-switches that transforms $H$ into $H^*$, where $NDL(H^*)=NDL(H)$ and neighbors of $w_H$ in $H^*$  are adjacent to vertices with the same degrees as neighbors of $w_H$ in $H$, but with highest $NDL$. 

 \bigskip

Observe that the degrees of the neighbors of $w$ and $w_H$ in $G^*$ and $H^*$, respectively, are the same. In addition,
$$NDL_{H^*}(w_H)=NDL_H(w_H)=NDL_G(w)=NDL_{G^*}(w)$$ 
Consider $G'=G^*-w$ and $H'=H^*-w_H$.  Then $NDL(G')=NDL(H')$. 
By the induction hypothesis applied to $G'$ and $H'$, there exists a sequence  of degree restricted 2-switches  that transforms $G'$ to $H'$. These degree restricted 2-switches do not involve $w$ and $w_H$, which have the same $NDL$ in $G^*$ and $H^*$, respectively. So applying this sequence of degree restricted 2-switches transforms $G^*$ to $H^*$. Finally, we can transform $G$ to $H$ by transforming $G$ to $G^*$, then $G^*$ to $H^*$, and then (in reverse order) $H$ to $H^*$. $\qed$
 
 \bigskip

\begin{lemma} \label{lemma-2} Let $G$ be a graph on $n$ vertices with diameter $2$. Then $$deg(v)=n-1-|N_2(v)|.$$
\end {lemma}

 \bigskip

\noindent {\bf Proof.} Since $G$ has diameter 2, every vertex is of distance 1 or 2 from every other vertex. So for each $v\in V(G)$, $$V(G)=\{v\}\cup N(v)\cup N_2(v),$$ where 
$|N(v)|\cap |N_2(v)|=\phi$. Since $deg(v)=|N(v)|$, $$n=1+deg(v)+|N_2(v)|.$$ Therefore $$deg(v)=n-1 -|N_2(v)|.$$
$\qed$

 \bigskip

\begin{figure}[h]
\centering
\includegraphics[width=1.5in]{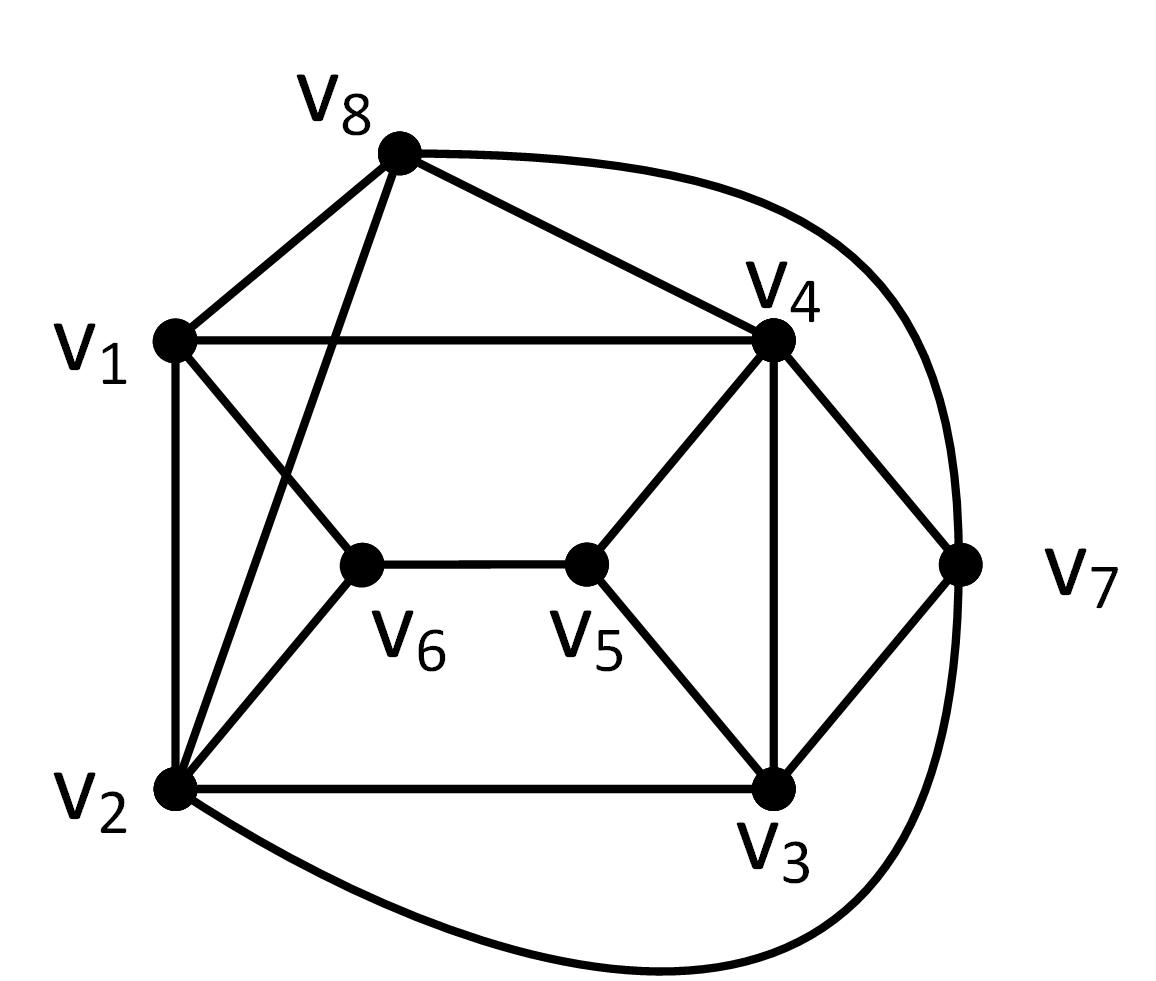}
\caption{ A diameter 2 graph
\label{fig04} }
\end{figure}

 The main idea in the proof of Theorem 1.1 is that if the graph has diameter 2, then  we can recover $NDL$ from $N_2DL$ and vice versa. Moreover, we can recover the degree sequence from $NDL$ in any graph. Let us look at an example to illustrate this point. Consider the  graph $G$ with diameter 2 shown in Figure \ref{fig04}. It has degree sequence
$5, 5, 4, 4, 4, 4, 3, 3$ and $NDL$ and $N_2DL$

\centerline{
\begin{tabular}{c|c|c}
&$NDL$&$N_2DL$\\
\hline
$v_1$ & $5, 5, 4, 3$    & $4, 4, 3$ \\
$v_2$ & $4, 4, 4, 4, 3$ & $5, 3$  \\
$v_3$ & $5, 5, 4, 3$    & $4, 4, 3$\\
$v_4$ & $4, 4, 4, 4, 3$ & $5, 3$\\
$v_5$ & $5, 4, 3$       & $5, 4, 4, 4$\\
$v_6$ & $5, 4, 3$       & $5, 4, 4, 4$\\
$v_7$ & $5, 5, 4, 4$    & $4, 3, 3$\\
$v_8$ & $5, 5, 4, 4$    & $4, 3, 3$
\end{tabular}}
\noindent Observe that $N_1(v_1)=\{v_2, v_4, v_6, v_8\}$. So the members of $N_2(v_1)$ are the rest of the vertices (except $v_1$ itself). Thus $N_2(v_1)=\{v_3, v_5, v_7\}$

 \bigskip

\noindent {\bf Proof of Theorem \ref{diam2}.} 
Suppose $G$ and $H$ have the same $N_2DL$. Then for every vertex $v_G$ in $G$, there is a vertex $v_H$ in $H$ with the same $N_2DL$ and vice-versa. Thus there is a one-to-one correspondence between the vertices of $G$ and $H$ such that for every pair of corresponding vertices $v_G$ and $v_H$, $$N_2DL(v_G)=\{deg_G (u) \ |\  u\in N_2(v_G)\}, $$ $$N_2DL(v_H)=\{deg_H (u) \  | \ u\in N_2(v_H)\},$$ and 
$N_2DL(v_G)=N_2DL(v_H)$. 
Observe  that $|N_2(v_G)|$ is the number of entries in $N_2DL(v_G)$ and $|N_2(v_H)|$ is the number of entries in $N_2DL(v_H)$. Therefore $|N_2(v_G)|=|N_2(v_H)|$. By Lemma \ref{lemma-2}  $$deg_G(v_G)=n-1 -|N_2(v_G)|=n-1-|N_2(v_H)|=deg_H(v_H).$$  Thus the degree sequence of $G$ and $H$ can be obtained from $N_2DL$. Moreover the degree sequence of $G$ and $H$ are the same. 

 \bigskip

Next observe that if $N_2DL(v) =\{deg(u) \ | \  u \in N_2(v)\}$, then since the degree sequence  is known, $NDL(v)=\{deg(u) \ | \  u \not \in \{v\}\cup N_2(v)\}$. Since $G$ and $H$ have the same $N_2DL$ and the same degree sequence, $G$ and $H$ must have the same $NDL$. Lemma \ref{NDL} implies that   $G$ can be transformed into $H$ by a sequence of degree restricted $2$-switches. 

 \bigskip
 
Conversely, suppose $G$ can be transformed into $H$ by a sequence of degree restricted $2$-switches. Lemma \ref{NDL} implies that $NDL$ is maintained at each stage (even if the diameter changes) so at the end of the sequence of degree restricted $2$-switches $G$ and $H$ have the same $NDL$.  Thus there is a one-to-one correspondence between the vertices of $G$ and $H$ such that for every pair of corresponding vertices $v_G$ and $v_H$, $$NDL(v_G)=\{deg_G (u) | u\in N(v_G)\},$$ $$NDL(v_H)=\{deg_H (u) | u\in N(v_H)\}$$ and 
$NDL(v_G)=NDL(v_H)$.  
Observe  that $|N(v_G)|$ is the number of entries in $NDL(v_G)$ and $|N(v_H)|$ is the number of entries in $NDL(v_H)$. Therefore $|N(v_G)|=|N(v_H)|$ and $deg(v_G)=deg(v_H)$. Thus the degree sequence of $G$ and $H$ can be obtained from $N_2DL$ and they are the same. 
 \bigskip

Next, observe that if $NDL(v) =\{deg(u) \ | \  u \in N(v)\}$, then since the degree sequence  is known, $N_2DL(v)=\{deg(u) \ | \  u \not \in \{v\}\cup N(v)\}$. In conclusion, if  $G$ and $H$ have the same $NDL$ and  the same degree sequence, then since $G$ and $H$ are diameter 2 graphs they must have the same $N_2DL$. $\qed$

 \bigskip


\noindent {\bf References}

 \bigskip

\begin{enumerate}

\item  \label{BarrusDonovan2018} M. Barrus and E. Donovan (2018). Neighborhood degree lists of graphs, {\it Discrete Mathematics}, {\bf 341(1)}, 175-183.

\item \label{Bassleretal2015} K.E. Bassler, C. I. Del Genio, P. L. Erdos, I. Miclos, Z Toroczkai (2015). Exact sampling of graphs with prescribed degree correlations, {\it New Journal of Physics} {\bf 17}.

\item \label{BuckleyHarary1990} F. Buckley and F. Harary (1990). Distance in graphs, Addison-Wesley.

\item \label{Godsil1995} C. D. Godsil (1995). Problems in Algebraic Combinatorics, {\it The Electronic Journal of Combinatorics}   2.

\item \label{FulkersonHoffmanMcAndrew1965} D. R. Fulkerson, A. J. Hoffman, and M.H. McAndrew (1965). Some properties of graphs with multiple edges, {\it Canadian Journal of Mathematics}, {\bf 17}, 166 - 177.

\item \label{Havel1955} V. Havel (1955).  A remark on the existence of finite graphs, {\it Cas. Pest. Mat.} {\bf 80}, 477–480

\item \label{Hakimi1962} S. L. Hakimi (1962).  On realizability of a set of integers as degrees of the vertices of a linear graph. I, {\it Journal of the Society for Industrial and Applied Mathematics}, {\bf 10}, 496–506.

\item \label{Randic1979} M. Randic (1979). Characterization of atoms, molecules, and classes of molecules based on paths enumerations. MATCH {\bf 7}, 5 - 64.

\item \label{Riberioetal2017} L. F. R. Robeiro, P. H. P. Saverese, and D. R Figueiredo (2017), Struc2vec: Learning Node Representations from Structural Identity, {\it Proceedings of the $23$rd ACM SIGKDD International Conference on Knowledge Discovery and Data Mining}  385-394.

\item \label{DamKoolen2016} E. R. van Dam and J. H. Koolen (2016) Distance-regular graphs,  {\it The Electronic Journal of Combinatorics} DS22.

\item \label{West1995} D. B. West (1995). {\it Introduction to Graph Theory}, Second Edition, Pearson.

\end{enumerate}

\end {document}